\documentclass{article}             
\usepackage{amssymb}

\def\bR{{\mathbb{R}}}

\def\fg{{\mathfrak{g}}}
\def\fh{{\mathfrak{h}}}
\def\fm{{\mathfrak{m}}}

\def\bA{{\mathbb{A}}}
\def\bB{{\mathbb{B}}}

\newtheorem       {theorem}{Theorem}
\newtheorem{defin}[theorem]{Definition}

\newtheorem{lemma}[theorem]{Lemma}

\newtheorem{con}   [theorem]{Conjecture}

\begin{document}

\begin{center}
{\Large {Structure of geodesics for Finsler metrics\\
arising from Riemannian g.o. metrics}}
\bigskip

{\large{Teresa Arias-Marco${}^1$ and Zden\v ek Du\v sek${}^2$} }
\bigskip

${}^1$ Department of Mathematics, University of Extremadura,\\
Av. de Elvas s/n, 06006 Badajoz, Spain\\
ariasmarco@unex.es\\
(ORCID: 0000-0003-0984-0367)

${}^2$ Institute of Technology and Business in \v Cesk\'e Bud\v ejovice\\
Okru\v zn\'\i\ 517/10, 370 01 \v Cesk\'e Bud\v ejovice, Czech Republic\\
zdusek@mail.vstecb.cz, corresponding author\\
(ORCID: 0000-0003-3073-9562)
\end {center}

\bigskip
 
\begin{abstract} 
Homogeneous geodesics of homogeneous Finsler metrics derived from two or more
Riemannian geodesic orbit metrics are investigated.
For a broad newly defined family of positively related Riemannian geodesic orbit metrics,
geodesic lemma is proved and
it is shown that the derived Finsler metrics have also geodesic orbit property.
These Finsler metrics belong to the newly defined class of the $\alpha_i$-type metrics
which includes in particular the $(\alpha_1,\alpha_2)$ metrics.
Geodesic graph for the sphere
${\mathrm{S}}^7={\mathrm{Sp(2)}}{\mathrm{U}}(1)/{\mathrm{Sp(1)}{\mathrm{diag}}{\mathrm{U}}(1)}$
with geodesic orbit Finsler metrics of the new type $(\alpha_1,\alpha_2,\alpha_3)$,
arising from two or more Riemannian geodesic orbit metrics,
is analyzed in detail.
This type of metrics on $S^7$ is one of the missing cases in a previously published classification
of geodesic orbit metrics on spheres.
\end{abstract}
\bigskip
 
\noindent
{\bf MSClassification:} {53C22, 53C60, 53C30}

\section{Introduction}
On a smooth manifold, Finsler metrics of the type $(\alpha_1,\alpha_2)$ were introduced
by S. Deng and M. Xu in \cite{DX} as metrics of the form
\begin{eqnarray}
\label{alpha}
F(y) & = & f(\alpha_1(y_1), \alpha_2(y_2)),
\end{eqnarray}
where $\alpha_i$ are nondegenerate symmetric bilinear forms on complementary subbundles of the tangent bundle
and $y_i$ are projections of vectors $y\in TM$ onto these subbundles.
These metrics were discussed also in some further papers,
but we did not find analysis of the function $f$
(conditions which $f$ must satisfy to obtain Finsler metric $F$),
neither nontrivial examples of such metrics.
On the other hand, in the paper \cite{JS}, M.A. Javaloyes and M. S\'anchez considered metrics of the type
\begin{eqnarray}
\label{ma}
F(y) & = & \sqrt{L(F_1(y),\dots, F_k(y),\beta_1(y),\dots,\beta_l(y))},
\end{eqnarray}
where $F_i$ are Finsler metrics and $\beta_j$ are one-forms on $M$.
They derived conditions for the function $L$ (to obtain the Finsler metric $F$)
and these conditions allow to construct many explicit examples of such metrics.
We are going to consider a special case of metrics (\ref{ma})
on homogeneous manifolds. They naturally provide examples of metrics of type (\ref{alpha})
and lead also to their generalization.

A homogeneous Finsler manifold $(M,F)$
admits a connected Lie group $G$ of isometries which acts transitively on $(M,F)$.
If we choose an origin $p\in M$ and denote by $H$ the isotropy group at $p$,
the manifold $(M,F)$ can be naturally identified with the homogeneous space $(G/H,F)$.
We denote by $\fg$ and $\fh$ be the Lie algebras of $G$ and $H$, respectively.
One can find a reductive decomposition $\fg=\fm+\fh$,
where $\fm\subset\fg$ is an ${\rm Ad}(H)$-invariant vector subspace
for the adjoint representation ${\mathrm{Ad}}\colon H\times\fg\rightarrow\fg$ of $H$ on $\fg$.
If we fix a reductive decomposition,
the vector space $\fm$ is naturally identified with the tangent space $T_pM$.
Using this identification, given the Finsler metric $F$ on $M$, we obtain the ${\mathrm{Ad}}(H)$-invariant
Minkowski norm and the ${\mathrm{Ad}}(H)$-invariant fundamental tensor on $\fm$.
Standard concepts in Finsler geometry can be found in monographs, for instance
\cite{BCS} by D. Bao, S.-S. Chern and Z. Shen,
\cite{SS} by Y.-B. Shen and Z. Shen,
or \cite{De} by S. Deng.
For the notion of derivative along a curve and geodesics, we use the Chern connection.

A geodesic $\gamma(t)$ through $p$ is {\it homogeneous} if it is an orbit
of a one-parameter subgroup of the group of isometries $G=I_0(M)$.
If $(M,F)$ is identified with the homogeneous space $(G/H,F)$,
there exists a nonzero {\it geodesic vector} $X\in\fg$ such that $\gamma(t)={\rm exp}(tX)(p)$ for all $t\in\bR$.
The homogeneous space $(G/H, F)$ is called a {\it geodesic orbit space} or just {\it g.o. space},
if every geodesic of $(G/H, F)$ is homogeneous.
Homogeneous geodesics and g.o. spaces were studied broadly in Riemannian geometry by many authors,
see for example \cite{AA}, \cite{AN}, \cite{ASS1}, \cite{CCZ}, \cite{DuS2}, \cite{DKN},
\cite{GoNi}, \cite{KNi}, \cite{Lu}, \cite{Ni} and the references therein,
or the monograph \cite{BN} by V. Berestovskii and Yu. Nikonorov.
One of the techniques for the study of g.o. spaces is based on the concept of geodesic graphs,
proposed by J. Szenthe in \cite{Sz}, in the affine setting.
This concept was further developed by O. Kowalski and S. Nik\v cevi\'c in \cite{KNi} for the Riemannian setting.
A {\it geodesic graph} is an ${\mathrm{Ad}}(H)$-equivariant map $\xi\colon\fm\rightarrow\fh$
such that $X +\xi(X)$ is a geodesic vector for each $o\neq X\in\fm$.
We are interested in the algebraic structure of the geodesic graph.
The existence of a linear geodesic graph is equivalent with the natural reductivity of the space $G/H$.
In general, the components of the Riemannian geodesic graph are rational functions $\xi_i=P_i/P$,
where $P_i$ and $P$ are homogeneous polynomials and ${\mathrm{deg}}(P_i)={\mathrm{deg}}(P)+1$.
See for example the recent survey paper \cite{DuS2} by the second author for
details and references to geodesic graphs on various classes of Riemannian g.o. manifolds.

Recently, g.o. spaces attained attention in Finsler geometry.
In \cite{YD}, Z. Yan and S. Deng studied relations of some Finsler g.o. spaces with Riemannian g.o. spaces,
for example for Randers g.o. metrics on spheres $S^{2n+1}={\mathrm{U}}(n+1)/{\mathrm{U}}(n)$
and for weakly symmetric metrics on some nilpotent groups.
In \cite{DuCMUC} and \cite{DuWS}, the second author investigated
Randers geodesic orbit metrics and special families of weakly symmetric Finsler metrics
on modified H-type groups using geodesic graphs.
Relations of Finslerian geodesic graphs with the Riemannian geodesic graphs were observed.
In \cite{DuNew}, second author studied general homogeneous Finsler $(\alpha,\beta)$ metrics,
which cover Randers metrics as a special case.
It was proved that all $(\alpha,\beta)$ metrics arising from a Riemannian g.o. metric $\alpha$
and an invariant one-form $\beta$ are Finsler g.o. metrics, possibly with respect to an extended
isometry group. See also the paper \cite{Y} by Z. Yan for a study of $(\alpha,\beta)$ metrics.
In the paper \cite{AMD}, present authors illustrated constructions of geodesic graphs
for geodesic orbit $(\alpha,\beta)$ metrics on spheres.
Classification list of these metrics on spheres and on projective spaces was given.
Geodesic orbit Finsler metrics on spheres were studied also by M. Xu in \cite{Xu}.
A classification list of types of these metrics for each possible isometry group was given,
with some missing cases.
In the present paper, we analyze the g.o. property of Finsler metrics of the type (\ref{ma}),
without one-forms and with the assumption that $F_i$ are norms determined by Riemannian g.o. metrics $g_i$,
hence $F_i=\sqrt{g_i}$.
We focus on the natural family of positively related initial Riemannian g.o. metrics, see Definition \ref{d2}.
We derive geodesic lemma for the characterization of geodesic vectors
and we prove that the Finsler metric $F$ given by the formula
\begin{eqnarray}
\label{ma2}
F(y) & = & \sqrt{L(\sqrt{g_1(y)},\dots, \sqrt{g_k(y)})},
\end{eqnarray}
is also a geodesic orbit metric.
It will be observed that these metrics provide natural examples of metrics of the type (\ref{alpha}).
This type can be further generalized to the new $\alpha_i$-type metrics, see Definition \ref{ai} below.
With an example of $S^7$, we illustrate the construction of geodesic graph in detail for metrics given by formula
(\ref{ma2}), which appear to be also of the type $(\alpha_1,\alpha_2,\alpha_3)$.
Metrics of this type are missing in the classification list in \cite{Xu}. 

\section{Preliminaries}
Let us consider a homogeneous manifold $M$ expressed as a homogeneous space $G/H$
with a fixed reductive decomposition $\fg=\fh+\fm$.
Let $p\in M$ and $F_p$ be an ${\mathrm{Ad}}(H)$-invariant Minkowski norm on $\fm\simeq T_pM$.
For any $q\in M$, we obtain a Minkowski norm $F_q$ at $T_qM$ by the formula
\begin{eqnarray}
\nonumber
F_q(\sigma_*X) = F_p(X), \qquad X\in T_pM,
\end{eqnarray}
where $\sigma\in G$ and $\sigma(p)=q$. It leads to the Finsler metric $F$ on $M$.
We shall work mostly just with the vector space $\fm$.
We denote the ${\mathrm{Ad}}(H)$-invariant
Minkowski norm on $\fm$ also by $F$ and its fundamental tensor on $\fm$ by $g$.
\begin{defin}
\label{ai}
Let $(G/H,F)$ be a homogeneous Finsler space with a reductive decomposition $\fg=\fm+\fh$.
Consider the adjoint action of $H$ on $\fm$, which leads to the irreducible decomposition
\begin{eqnarray}
\nonumber
\fm & = & \oplus_{i=1}^s \fm_i.
\end{eqnarray}
Choose symmetric positively definite ${\mathrm{Ad}}(H)$-invariant bilinear forms $\alpha_i$ on $\fm_i, i=1\dots s$
and let $y_i$ be the corresponding projections of a vector $y\in\fm$ onto $\fm_i$.
The Minkowski norm $F$ on $\fm$ and the corresponding homogeneous Finsler metric on $G/H$ is of the $\alpha_i$-type
if there exist a smooth function $f\colon [ 0,\infty)^s\rightarrow{\mathbb{R}}$ such that
\begin{eqnarray}
\label{f5}
F(y) & = & f(\alpha_1(y_1),\dots, \alpha_s(y_s)),\qquad y\in\fm.
\end{eqnarray}
\end{defin}
This definition obviously covers the homogeneous $(\alpha_1,\alpha_2)$ metrics from \cite{DX}.
However, as far as we know, necessary and sufficient conditions for the function $f$,
which imply that the function $F$ is a Minkowski norm, were never investigated.
We are also not aware of any explicit and nontrivial examples of these metrics.

Another approach how to construct new Finsler metrics from given Finsler metrics was proposed in \cite{JS},
see formula (\ref{ma}) above.
For the simplicity and technical reasons later, we shall avoid the one-forms and
consider given metrics $F_i$ to be Riemannian. Hence we write the new Minkowski norm $F$ on $\fm$ in the form
(\ref{ma2}), where $g_i$ are scalar products on $\fm$.
According to Theorem 4.1 in \cite{JS}, to obtain a Minkowski norm $F$, the continuous function $L$
must satisfy:\\
(i) be smooth and positive away from $0$,\\
(ii) be positively homogeneous of degree 2,\\
(iii) $L_{,i} \geq 0$, for $i=1,\dots,k$,\\
(iv) $Hess(L)$ be positive semi-definite,\\
(v) $L_{,1} + \dots + L_{,k} >  0$.\\
Here the comma in $L_{,i}$ means the derivative with respect to the corresponding coordinate.
Obviously, there are many functions which satisfy these conditions,
for example, one can choose $L=\sqrt{g_1}+\dots+\sqrt{g_k}$.
\begin{defin}
\label{d2}
Let $G/H$ be a homogeneous space with a reductive decomposition $\fg=\fm+\fh$.
Consider the ${\mathrm{Ad}}(H)$-invariant irreducible decomposition $\fm=\oplus_{i=1}^s \fm_i$
and let $\alpha_i$ be ${\mathrm{Ad}}(H)$-invariant scalar products on the respective spaces $\fm_i$.
We consider the family of scalar products
\begin{eqnarray}
\nonumber
g(c_1,\dots,c_s) & = & \sum_{i=1}^s c_i\cdot \alpha_i,
\end{eqnarray}
for any numbers $0<c_i\in {\mathbb{R}}$.
This family of scalar products on $\fm$ and corresponding family of Riemannian metrics on $G/H$
will be called scalar products positively related and metrics positively related.
\end{defin}
If the scalar products $g_j, j=1,\dots,k$, are positively related,
then obviously the Finsler metric $F$ given by the formula (\ref{ma2}) can be written also using formula (\ref{f5}),
for some function $f$ determined by $L$, and hence it is of the $\alpha_i$-type.
We are going to analyze geodesic orbit property for this type of Minkowski norms on $\fm$
and corresponding Finsler metrics on $M$.
Throughout the rest of the paper, we always assume that the function $L$ satisfies the conditions
(i)-(v) above.

Let us recall the geodesic lemma, which was proved for general Finsler metrics by D. Latifi in \cite{La}.
\begin{lemma}%[\cite{La}]
\label{golema2}
Let $(G/H,F)$ be a homogeneous Finsler space with a reductive decomposition $\fg=\fm+\fh$.
A nonzero vector $y\in{\fg}$ is geodesic vector if and only if it holds
\begin{equation}
\label{gl2}
g_{y_\fm} ( y_{\mathfrak m}, [y,u]_{\mathfrak m} ) = 0 \qquad \forall u\in{\mathfrak m},
\end{equation}
where the subscript $\fm$ indicates the projection of a vector from $\fg$ to $\fm$.
\end{lemma}

\section{Main results}
We shall first adapt formula (\ref{gl2})
in terms of the Riemannian metrics $g_i$ and the function $L$.
First, we need the expression of the fundamental tensor $g$.
\begin{lemma}
\label{lem1}
Let $g_1,\dots, g_k$ be homogeneous Riemannian metrics on $G/H$
and let $\fg=\fm+\fh$ be a reductive decomposition.
We use the same notation for the bilinear forms on $\fm$.
Let $F=\sqrt{L(\sqrt{g_1},\dots,\sqrt{g_k})}$ on $\fm$, which gives a homogeneous Finsler metric on $G/H$.
For arbitrary vectors $y,v\in\fm$, the fundamental tensor $g$ of $F$ satisfies the formula
\begin{eqnarray}
\label{fg}
g_y(y,v) & = & \sum_{j=1}^k B_j(y) \cdot g_j(y,v),
\end{eqnarray}
where the functions $B_j(y)$ are given by
\begin{eqnarray}
\label{fb}
B_j(y) & = & \frac{ L_{,j}}{2\,\sqrt{g_j(y,y)}} = \frac{ L_{,j}}{2\,F_j(y)}, \quad j=1\dots k.
\end{eqnarray}
\end{lemma}
{\it Proof.}
Let $y,v\in\fm$ be arbitrary fixed vectors.
For arbitrary $t\in\bR$, for the vector $w=y+tv\in\fm$, it holds
\begin{eqnarray}
\nonumber
F^2(w) & = & L(\sqrt{g_1(w,w)},\dots, \sqrt{g_k(w,w)}).
\end{eqnarray}
We now differentiate with respect to $t$ and apply the chain rule.
We obtain
\begin{eqnarray}
\label{F1}
\nonumber
\frac{\partial F^2(w)}{\partial t} & = & \sum_{j=1}^k L_{,j} \frac{1}{\sqrt{g_j(w,w)}}g_j(w,v)
\end{eqnarray}
and
\begin{eqnarray}
\nonumber
2\,g_y(y,v) & = &
\frac{\partial F^2(w)}{\partial t}\Big |_{t=0} =
\sum_{j=1}^k L_{,j} \frac{1}{\sqrt{g_j(y,y)}} g_j(y,v).
\end{eqnarray}
Using the notation (\ref{fb}), we obtain the statement.
$\hfill\square$

Let us now consider the case when the metrics $g_j, j=1\dots k,$
are positively related Riemannian g.o. metrics.
Let $\fg=\fm+\fh$ be a reductive decomposition with the ${\mathrm{Ad}}(H)$-invariant irreducible decomposition
$\fm=\oplus_{i=1}^s \fm_i$.
According to Definition \ref{d2}, scalar products $g_j$ have corresponding decompositions
\begin{eqnarray}
\label{f4}
g_j & = & \sum_{i=1}^s a_{ji}\cdot \alpha_i, \qquad j=1\dots k,
\end{eqnarray}
where $0<a_{ji}\in{\mathbb{R}}$ and
$\alpha_i$ are some initial ${\mathrm{Ad}}(H)$-invariant scalar products on $\fm_i$.
Our experience with g.o. metrics leads us to the following conjecture, 
whose proof in full generality seems to be not easy.
\begin{con}
\label{con1}
Consider a family of positively related Riemannian metrics.
If one metric of this family is a g.o. metric, then all metrics from this family are also g.o. metrics.
\end{con}
We remark that the families of g.o. metrics may be broader than a family of positively related metrics.
In some cases, in individual scalar products $\alpha_i$, symmetric ${\mathrm{Ad}}(H)$-invariant
operators may be applied. See for example \cite{Lu} for the case of H-type groups
or \cite{AA} or \cite{AN} for compact spaces.
But we keep further constructions simpler by considering only families of positively
related g.o. metrics, as above.

Denote now
\begin{eqnarray}
\label{fc}
C_i(y) & = & \sum_{j=1}^k{ B_{j}(y)} \cdot a_{ji}.
\end{eqnarray}
From formula (\ref{fg}) and (\ref{f4}), we obtain that the fundamental
tensor of the homogeneous Finsler metric $F=\sqrt{L(\sqrt{g_1},\dots,\sqrt{g_k})}$ satisfies
\begin{eqnarray}
\label{11}
g_y(y,v) & = & \sum_{j=1}^k B_j(y) \cdot g_j(y,v)  \cr
& = & \sum_{j=1}^k B_j(y) \cdot \sum_{i=1}^s a_{ji} \cdot \alpha_i(y,v) =
\sum_{i=1}^s C_i(y) \cdot \alpha_i(y,v).
\end{eqnarray}
\begin{lemma}
\label{c2}
Let $G/H$ be a homogeneous space with a reductive decomposition $\fg=\fm+\fh$
and with the ${\mathrm{Ad}}(H)$-irreducible decomposition $\fm=\oplus_{i=1}^s \fm_i$.
Let $g_j$ be positively related Riemannian metrics 
with decompositions $(\ref{f4})$.
Let $F=\sqrt{L(\sqrt{g_1},\sqrt{g_2},\dots,\sqrt{g_k})}$
with the corresponding decomposition $(\ref{fg})$ of the fundamental tensor $g$
and determine functions $C_i(y)$ by formula $(\ref{fc})$.
The vector $y+\xi(y)$, where $y\in\fm$ and $\xi(y)\in\fh$,
is geodesic vector for the Finsler metric $F$ if and only if
\begin{eqnarray}
\label{c7}
\sum_{i=1}^s C_i(y) \cdot \alpha_i\Bigl (y,[y+\xi(y),U]_{\mathfrak m}\Bigr ) & = & 0, \quad \forall U\in\fm.
\end{eqnarray}
\end{lemma}
{\it Proof.}
Follows immediately from Lemma \ref{golema2}, Lemma \ref{lem1} and notations above.
$\hfill\square$

\begin{theorem}
Let $G/H$ be a homogeneous space with a reductive decomposition $\fg=\fh+\fm$ and with the
${\mathrm{Ad}}(H)$-irreducible decomposition $\fm=\oplus_{i=1}^s \fm_i$.
Let $g_j, j=1\dots k,$ be from a family of positively related Riemannian metrics on $G/H$,
all of which are g.o. metrics.
Then any homogeneous Finsler metric of the type
$F=\sqrt{L(\sqrt{g_1},\sqrt{g_2},\dots,\sqrt{g_k})}$ on $G/H$
is also a g.o. metric.
\end{theorem}
{\it Proof.}
We consider condition (\ref{c7}) and we fix the arbitrary vector $y\in\fm$.
The values $C_i(y)$ are positive real numbers and the Riemannian metric
\begin{eqnarray}
\nonumber
C_1(y) \cdot \alpha_1 + \dots + C_s(y) \cdot \alpha_s
\end{eqnarray}
belongs to the family of initial positively related g.o. metrics.
It follows that, in the condition (\ref{c7}), there exist a vector $\xi(y)$ such that (\ref{c7})
is satisfied.
Because $y\in\fm$ was arbitrary, the metric $F$ is also a g.o. metric.
$\hfill\square$

\section{Example
$S^{7}={\mathrm{Sp}}(2)\cdot{\mathrm{U}}(1)/{\mathrm{Sp}}(1)\cdot{\mathrm{diag}}({\mathrm{U}}(1))$ }
On the Lie algebra level, we choose a basis $\{H_1,H_2,H_3\}$ of the Lie algebra $\fh={\mathfrak{sp}}(1)$,
such that, with respect to this basis, an element of $\fh$ is given by a matrix
\begin{eqnarray}
\nonumber
\left (
\begin{array}
{cccc}
ih_1 & -h_2-ih_3 & 0 & 0 
\\{\medskip}
h_2-ih_3 & -ih_1 & 0 & 0
\\{\medskip}
0 & 0 & 0 & 0 
\\{\medskip}
0 & 0 & 0 & 0 
\end{array}
\right ).
\end{eqnarray}
The reductive complement $\fm$ in the decomposition ${\mathfrak{sp}}(2)= {\mathfrak{sp}}(1) + \fm$
and its basis $B=\{X_1,\dots,X_4,Z_1,Z_2,Z_3\}$ will be chosen such that,
with respect to this basis $B$, each element of $\fm$ is given by the matrix
\begin{eqnarray}
\nonumber
\left (
\begin{array}
{cccc}
0 & 0 & x_1 + ix_2 & -x_3 - ix_4
\\{\medskip}
0 & 0 & x_3 - ix_4 & x_1 - ix_2
\\{\medskip}
- x_1 + ix_2 & -x_3 - ix_4 & iz_1 & -z_2 - ix_3
\\{\medskip}
x_3 - ix_4 & -x_1 - ix_2 & z_2 - iz_3 & -iz_1
\end{array}
\right ).
\end{eqnarray}
It is straightforward to determine the Lie brackets on $\fg\simeq{\mathfrak{sp}}(2)$.
One can find explicit formulas for example in \cite{AMD},
where geodesic graphs for $(\alpha,\beta)$ metrics on this space were calculated,
in a modified reductive decomposition.
The adjoint action of $\fh$ on $\fm$ is given by the operators
\begin{eqnarray}
\nonumber
{\rm ad}(H_1)|_\fm = A_{12}+A_{34},\quad
{\rm ad}(H_2)|_\fm = A_{13}-A_{24},\quad
{\rm ad}(H_3)|_\fm = A_{14}+A_{23}.
\end{eqnarray}
We further introduce the operator
\begin{eqnarray}
\nonumber
& W={\rm ad}(Z_1)|_\fm = 2B_{23}-A_{12}+A_{34}.
\end{eqnarray}
Here $A_{ij}$ and $B_{ij}$ are the elementary operators acting nontrivially by
$A_{ij}(X_i)=X_j$ and $A_{ij}(X_j)=-X_i$, resp.
$B_{ij}(Z_i)=Z_j$ and $B_{ij}(Z_j)=-Z_i$.
We put $\widetilde\fh={\mathrm{span}}\{\fh,W\}\simeq{\mathrm{sp(1)}}\oplus {\mathrm{u(1)}}$.

We will now consider the family of invariant scalar products with respect to $\widetilde\fh$
and corresponding invariant Riemannian metrics.
Obviously, the ${\mathrm{Ad}}(\widetilde H)$-invariant subspaces of $\fm$ are
$\fm_1={\mathrm{span}}(X_1,\dots,X_4)$, $\fm_2={\mathrm{span}}(Z_1)$ and $\fm_3={\mathrm{span}}(Z_2,Z_3)$.
If we choose any ${\mathrm{Ad}}(\widetilde H)$-invariant scalar products $\alpha_i$ on $\fm_i$, then any scalar product
$g=\sum_{i=1}^3 c_i\alpha_i$ for $0<c_i\in {\mathbb{R}}$ is ${\mathrm{Ad}}(\widetilde H)$-invariant.
It is well known that all these scalar products determine invariant Riemannian g.o. metrics on $\widetilde G/\widetilde H$
and we can observe this fact also from the further considerations.
Clearly, metrics from this 3-parameter family are positively related.

Let us now use the basis $B$ above and the natural scalar product $g_0$ in which this basis is orthonormal.
We put $\alpha_i=g_0|_{\fm_i}$.
Let the scalar products $g_j, j=1\dots k$, be $g_j=\sum_{i=1}^3 a_{ij}\alpha_i$ for some fixed $0<a_{ij}\in {\mathbb{R}}$.
Consider a Finsler metric $F=\sqrt{L(\sqrt{g_1},\sqrt{g_2},\dots,\sqrt{g_k})}$.
Note that $F$ can be alternatively written in the form $F=f(\alpha_1,\alpha_2,\alpha_3)$
for some function $f$ determined by $L$.
Hence the Finsler metric $F$ belongs to the family of the $\alpha_i$-type metrics.
In this case, it is of the $(\alpha_1, \alpha_2, \alpha_3)$ type.

For the fundamental tensor of $F$, we have
\begin{eqnarray}
\nonumber
g_y(y,v) & = & \sum_{j=1}^k B_j(y) \cdot g_j(y,v) \cr
         & = & C_1(y) \cdot \alpha_1(y,v) + C_2(y) \cdot \alpha_2(y,v) + C_3(y) \cdot \alpha_3(y,v),
\end{eqnarray}
according to formula (\ref{11}). We use formula (\ref{c7}) in the form
\begin{eqnarray}
\nonumber
\sum_{i=1}^3 C_i(y) \cdot \alpha_i\Bigl (y,[\xi(y),u]_{\mathfrak m}\Bigr ) & = &
- \sum_{i=1}^3 C_i(y) \cdot \alpha_i\Bigl (y,[y,u]_{\mathfrak m}\Bigr ) , \quad \forall u\in\fm.
\end{eqnarray}

\noindent
We write each vector $y\in\fm$ and each vector $\xi(y)\in\widetilde\fh$ in the forms
\begin{eqnarray}
\nonumber
y & = & x_1X_1 + \dots + x_4 X_4 + z_1 Z_1 + z_2 Z_2 + z_3 Z_3, \cr
\xi(y) & = & \xi_1(y) H_1 + \xi_2(y) H_2 + \xi_3(y) H_3 + \xi_4(y) W.
\end{eqnarray}
We obtain the system of linear equations for components $\xi_i=\xi_i(y)$, whose extended matrix
$(\bA | \bB)$ is
\begin{eqnarray}
\nonumber
\left (
\begin{array}
{cccc|c}
 {x_2}&{x_3}&{x_4}& {-x_2}&
 (1-2\frac{C_2}{C_1})z_1x_2 + (1-2\frac{C_3}{C_1})(z_2x_3 + z_3x_4)
\\{\medskip}
{-x_1}&{-x_4}&{ x_3}& { x_1}&
-(1-2\frac{C_2}{C_1})z_1x_1 + (1-2\frac{C_3}{C_1})(z_2x_4 - z_3x_3)
\\{\medskip}
{ x_4}&{-x_1}&{-x_2}& { x_4}&
-(1-2\frac{C_2}{C_1})z_1x_4 + (1-2\frac{C_3}{C_1})(- z_2x_1 + z_3x_2)
\\{\medskip}
 {-x_3}&{ x_2}&{-x_1}& {-x_3}&
 (1-2\frac{C_2}{C_1})z_1x_3 - (1-2\frac{C_3}{C_1})(z_2x_2 + z_3x_1)
\\{\medskip}
0 & 0 & 0 &  2z_3 & 2z_1z_3(\frac{C_2}{C_3}-1) 
\\{\medskip}
0 & 0 & 0 & -2z_2 & 2z_1z_2(1-\frac{C_2}{C_3})
\end{array}
\right ),
\end{eqnarray}
using the short notation $C_i=C_i(y)$. We denote further
\begin{eqnarray}
\nonumber
K_1=K_1(y)  = \frac{C_2}{C_3}-2\frac{C_2}{C_1}, \quad
K_2=K_2(y)  = 1-2\frac{C_3}{C_1}, \quad
K_3=K_3(y)  = \frac{C_2}{C_3}-1.
\end{eqnarray}
The unique solution of the above system of equations is
\begin{eqnarray}
\nonumber
\xi_1 & = &
\frac{
K_1 z_1 (x_1^2+ x_2^2- x_3^2- x_4^2 )+
2K_2 \bigl[
z_2 (x_2x_3 - x_1x_4) +
z_3 (x_1x_3 + x_2x_4) 
\bigr ]
}{ x_1^2+ x_2^2+ x_3^2+ x_4^2 }, \cr
\xi_2 & = &
\frac{
2K_1 z_1 (x_2x_3 + x_1x_4) +
K_2 \bigl[
z_2 (x_1^2- x_2^2+ x_3^2- x_4^2 ) +
2z_3 (x_3x_4 - x_1x_2)
\bigr ]
}{ x_1^2+ x_2^2+ x_3^2+ x_4^2 }, \cr
\xi_3 & = &
\frac{
2K_1 z_1 (x_2x_4 - x_1x_3) +
K_2 \bigl[
2z_2 (x_1x_2 + x_3x_4) +
z_3 (x_1^2- x_2^2- x_3^2+ x_4^2 )
\bigr ]
}{ x_1^2+ x_2^2+ x_3^2+ x_4^2 }, \cr
\xi_4 & = &  K_3 z_1.
\end{eqnarray}
We observe that, if we put $C_i(y)=c_i>0$, we obtain formulas for the geodesic graph
of the Riemannian metric $g(c_1,c_2,c_3)=\sum_{i=1}^3 c_i\alpha_i$.
The geodesic graph for the Finsler metric $F$ can be obtained
from the geodesic graphs for the family of positively related Riemannian metrics $g(c_1,c_2,c_3)$
by changing constants $c_i$ with the functions $C_i(y)$, according to formulas (\ref{fb}) and (\ref{fc}).

\section*{Statements and Declarations}
\textbf{Author contributions:} Both authors contributed equally
to this research and in writing the paper.

\noindent
\textbf{Funding:}
The research is supported by grant PID2019-10519GA-C22 funded by
Agencia Estatal de Investigaci\'on Espa\~{n}ola ($AEI/10.13$ $039/501100011033$).
The first author is also partially supported by grant GR21055 funded by Junta de Extremadura
and Fondo Europeo de Desarrollo Regional.

\noindent
\textbf{Conflicts of Interest:}
The authors have no financial or competing interests to declare that are relevant to the content of this article.


\begin{thebibliography}{[0000]}
\bibitem {AA} % ???
Alekseevsky, D., Arvanitoyeorgos, A.:
Riemannian flag manifolds with homogeneous geodesics,
{\it Trans. Am. Math. Soc.} {\bf{359,8}} (2007), 3769--3789.
\bibitem {AMD}
Arias-Marco, T., Du\v sek, Z.:
Geodesic graphs for geodesic orbit Finsler $(\alpha,\beta)$ metrics on spheres,
{\it Vietnam J. Math.}, to appear.
\bibitem {AN} % ???
Alekseevsky, D., Nikonorov, Yu.G.:
Compact Riemannian Manifolds with Homogeneous Geodesics,
{\it SIGMA} {\bf{5}} (2009), 093, 16pp.
\bibitem {ASS1} % ???
Arvanitoyeorgos, A., Souris, N.P., Statha, M.:
Geodesic orbit metrics in a class of homogeneous bundles over real and complex Stiefel manifolds,
{\it Geom. Dedic.} {\bf{215}} (2021), 31--50.
\bibitem {BCS}
Bao, D., Chern, S.-S., Shen, Z.:
An Introduction to Riemann-Finsler Geometry,
Springer Science+Business Media, New York, 2000.
\bibitem {BN}
Berestovskii, V., Nikonorov, Yu.:
Riemannian Manifolds and Homogeneous Geodesics,
Springer Nature Switzerland, Cham, 2020.
\bibitem {CCZ} % ???
Chen, H.B., Chen, Z.Q., Zhu, F.H.:
Geodesic orbit metrics on homogeneous spaces constructed by strongly isotropy irreducible spaces,
{\it Sci. China-Math.} {\bf{64}} (2021), 2313--2326.
\bibitem {De}
Deng, S.:
Homogeneous Finsler Spaces, Springer Science+Business Media, New York, 2012.
\bibitem {DX}
Deng, S., Xu, M.:
$(\alpha_1,\alpha_2)$-Metrics and Clifforf-Wolf Homogeneity,
{\it J. Geom. Anal.} {\bf{26}} (2016), 2282--2321.
\bibitem{DuS2}
Du\v sek, Z.:
Homogeneous geodesics and g.o. manifolds,
{\it Note Mat.} {\bf{38}} (2018), 1--15.
\bibitem{DuCMUC}
Du\v sek, Z.:
Geodesic graphs in Randers g.o. spaces,
{\it Comment. Math. Univ. Carol.} {\bf{61,2}} (2020), 195--211.
\bibitem{DuWS}
Du\v sek, Z.:
Structure of geodesics in weakly symmetric Finsler metrics on H-type groups,
{\it Arch. Math.-Brno} {\bf{56}} (2020), 265--275.
\bibitem{DuNew}
Du\v sek, Z.:
Geodesic orbit Finsler $(\alpha,\beta)$ metrics,
{\it Eur. J. Math.} {\bf{9,1}} (2023), 9.
\bibitem {DKN}
Du\v sek, Z., Kowalski, O. Nik\v cevi\'c, S.:
New examples of Riemannian g.o. manifolds in dimension $7$,
{\it Differ. Geom. Appl.} {\textbf{21}} (2004), 65--78.
\bibitem {GoNi}
Gordon, C.S., Nikonorov, Yu.G.:
Geodesic orbit Riemannian structures on ${\mathbb{R}}^n$,
{\it J. Geom. Phys.} {\bf{134}} (2018), 235--243.
\bibitem {JS}
Javaloyes, M.A., S\'anchez, M.:
On the definition and examples of Finsler metrics,
{\it Ann. Scuola Norm. Super. Pisa-Cl. Sci.} Vol. XIII (2014), 813--858.
\bibitem {KNi}
Kowalski, O., Nik\v cevi\' c, S.:
On geodesic graphs of Riemannian g.o. spaces,
{\it Arch. Math.} {\bf{73}} (1999), 223--234;
Appendix:
{\it Arch. Math.} \textbf{79} (2002), 158--160.
\bibitem {La}
Latifi, D.:
Homogeneous geodesics in homogeneous Finsler spaces,
{\it J. Geom. Phys.} {\bf{57}} (2007), 1421--1433.
\bibitem {Lu}
Lauret, J.:
Modified H-type groups and symmetric-like Riemannian spaces,
{\it Differ. Geom. Appl.} {\bf{10}} (1999), 121--143.
\bibitem {Ni}
Nikonorov, Yu.G.:
On the structure of geodesic orbit Riemannian spaces,
{\it Ann. Glob. Anal. Geom.} {\bf{52}}(3) (2017), 289--311.
\bibitem {SS}
Shen, Y.-B., Shen, Z.:
Introduction to Modern Finsler Geometry,
Higher Education Press, Beijing and World Scientific, Singapore, 2016.
\bibitem {Sz} Szenthe, J.:
Sur la connection naturelle \`a torsion nulle,
{\it Acta Sci. Math. (Szeged)} {\bf{38}} (1976), 383--398.
\bibitem {Xu}
Xu, M.:
Geodesic orbit spheres and constant curvature in Finsler geometry,
{\it Differ. Geom. Appl.} {\bf{61}} (2018), 197--206.
\bibitem {Y}
Yan, Z.:
Some Finsler spaces with homogeneous geodesics,
{\it Math. Nachr.} {\bf{290}},2-3 (2017), 474--481.
\bibitem {YD}
Yan, Z., Deng, S.:
Finsler spaces whose geodesics are orbits,
{\it Differ. Geom. Appl.} {\bf{36}} (2014), 1--23.
\end{thebibliography}
\end{document}